\pdfoutput=1
\documentclass[11pt,oneside]{amsart}
\usepackage[leqno]{amsmath}
\usepackage{amssymb,amsthm} 
\usepackage[parfill]{parskip} 
\usepackage{libertine} 
\usepackage[T1]{fontenc}
\usepackage{textcomp}
\usepackage[varqu,varl]{inconsolata} 
\usepackage{amsmath,amsthm}
\usepackage[libertine,bigdelims,vvarbb]{newtxmath}
\usepackage[supstfm=libertinesups,%
  supscaled=1.2,%
  raised=-.13em]{superiors}
\usepackage{graphicx} 
\usepackage{calc}

\newcommand{\beq}{\begin{equation}}
\newcommand{\eeq}{\end{equation}}
\numberwithin{equation}{section} 
\newtheorem{thm}[equation]{Theorem}
\newtheorem{lem}[equation]{Lemma}

\theoremstyle{remark}

\theoremstyle{definition}
\newtheorem{defn}[equation]{Definition}

\title{Lattice Exit Models}
\usepackage[paperwidth=7in, paperheight=10in, textwidth=4.75in, textheight=8in, includehead=true]{geometry}
\usepackage{url}
\usepackage[pdfborderstyle={/S/U/W 1}]{hyperref}
\def\bF{\mathbb{F}}

\def\bM{{\mathbb{M}}}

\def\bE{{\mathbb{E}}}
\def\setmax{\operatorname{setmax}}
\def\diag{\operatorname{diag}}

\let\emptyset\emptysetAlt
\hypersetup{
pdftitle={Terminal Search Multiverse},
pdfauthor={\textcopyright\ S. G. Williamson}
}
\author{S. Gill Williamson}
\thanks{Department of Computer Science and Engineering, 
University of California San Diego; \url{http://cseweb.ucsd.edu/~gill/}.
}

\date{05/01/2017}                                           
\begin{document}
\maketitle
\begin{abstract}  
We discuss a class of problems which we call {\em lattice exit models}.  At one level, these problems provide undergraduate level exercises in labeling the vertices of graphs (e.g., depth first search).  At another level (theorems about large scale regularities of  labels) they provide concrete geometric examples of ZFC independence. We note some combinatorial and algorithmic implications.
\end{abstract}

\section{Introduction}
We expand on some of the geometric and combinatorial concepts contained in the foundational work of H. Friedman (~\cite{hf:alc}, \cite{hf:nlc}).

Let $N$ be the set of nonnegative integers and $k\geq 2$.
For $z=(n_1, \ldots, n_k)\in N^k$, $\max\{n_i\mid i=1,\ldots, k\}$ will be
denoted by $\max(z)$.  Define $\min(z)$ similarly. 
\begin{defn}[Downward directed graph]
\label{def:dwndrctgrph}
Let $G=(N^k,\Theta)$ (vertex set $N^k$, edge set $\Theta$) be a directed graph.
If every $(x,y)$ of $\Theta$ satisfies $\max(x) >\max(y)$  then we call $G$ a {\em downward directed lattice graph}. 
\end{defn}
{\bf All lattice graphs that we consider will be {\em downward directed}.}\\
\begin{figure}[h]
\begin{center}
\includegraphics[scale=.85]{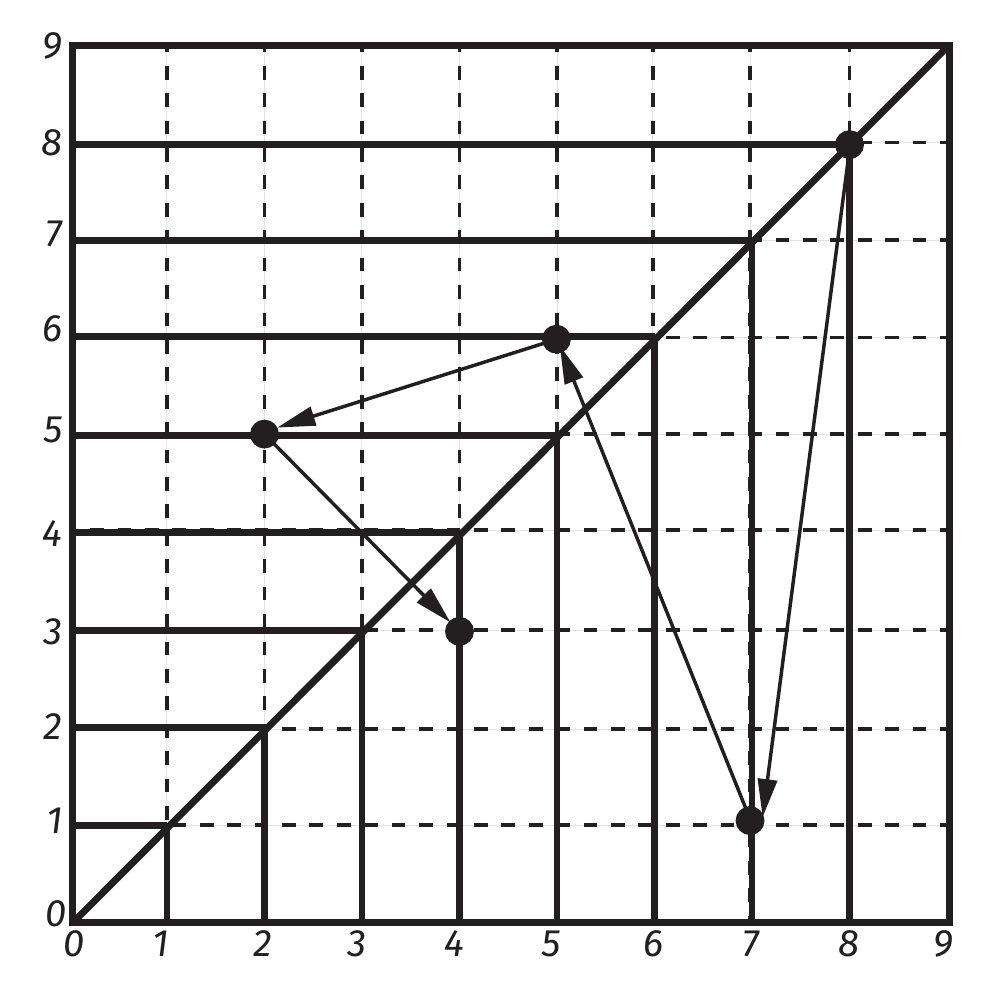}
\caption{Downward directed, $k=2$}
\label{fig:downdir}
\end{center}
\end{figure}
\begin{defn} [\bf Vertex induced subgraph $G_D$]
For $D\subset N^k$  let 
$G_D = (D, \Theta_D)$ be the subgraph of $G$ with vertex set $D$ and edge set 
$\Theta_D =\{(x,y)\mid (x,y)\in \Theta,\, x, y \in D\}$. We call $G_D$ the {\em subgraph of $G$ induced by $D$}. 
\label{def:vertexinduced}
\end{defn}
\vskip .25in  
\begin{defn} [\bf Path and terminal path in $G_D$]
\label{def:terminalpath}
For $t\geq 2$, a sequence of distinct vertices of $G_D$, $(x_1, x_2, \ldots, x_t)$, is a {\em path}
of length $t$ in $G_D$  if $(x_i, x_{i+1}) \in \Theta_D,\, i=1, \ldots\,, t-1$.  For $
x\in D$, $(x)$ as a path of length $1$. The path $(x_1, x_2, \ldots, x_t)$ is {\em terminal}
if there is no path of the form  $(x_1, x_2, \ldots, x_t, x_{t+1})$. 
We say $x$ is a {\em terminal vertex} of $G_D$ if the path $(x)$ is a terminal path in $G_D$.
\label{def:downpath}
\end{defn}

\begin{defn}[\bf Cubes and Cartesian powers in $N^k$]
The set  $E_1\times\cdots\times E_k$, where $E_i\subset N$, $|E_i|=p$, $i=1,\ldots, k,$   are $k$-cubes of length $p$.  If $E_i = E, i=1,\ldots, k,$ then this cube is  $E^k \coloneq \times^k E$, the $k$th Cartesian power of $E$.
\label{def:cubespowers}
\end{defn}

\begin{defn} [\bf $p_D$ and significant labels]
\label{def:pathlabel}
For  finite $D\subset N^k$, let $G_D = (D, \Theta_D)$.  Let $P_D(z)$ be the set of all $x$ with a path in $G_D$ from 
$z$ to $x$ (including the path $(z)$).
Define  ${p}_D$ by
$${p}_D(z) = \min(\{\min(x)\mid x\in P_D(z)\}).$$
We call ${p}_D$ the {\em total path label function}.
The set $\{z\mid {p}_D(z)<\min(z)\}$
is the {\em set of vertices with significant labels}.
The set $\{{p}_D(z)\mid {p}_D(z)<\min(z)\}$ is the {\em set of
significant labels} for ${p}_D.$ 
\label{def:alllabel}
\end{defn}

\vskip .25in
\begin{defn} [\bf $\hat{t}_D$ and significant labels]
\label{def:termlabel}
For finite $D\subset N^k$, let $G_D = (D, \Theta_D)$.  Let $T_D(z)$ be the set of all last vertices of terminal paths $(x_1, x_2, \ldots, x_t)$
where $z=x_1$.
Define $\hat{t}_D$ by
$$\hat{t}_D(z) = \max(z)\;\mathrm{if}\;(z)\;\mathrm{terminal, else}$$ 
$$\hat{t}_D(z)  = \min(\{\min(x)\mid x\in T_D(z)\})$$
We call $\hat{t}_D$ the 
{\em terminal path label function}.
The set $\{z\mid \hat{t}_D(z)<\min(z)\}$ is the set of vertices with significant labels.
The set $\{\hat{t}_D(z)\mid \hat{t}_D(z)<\min(z)\}$ is the set of
{\em significant labels} for $\hat{t}_D$.

\end{defn}

In definition~\ref{def:recursivedefs} we define a variation on $p_D$ of definition~\ref{def:pathlabel}
called $\hat{p}_D$ where
 $\hat{p}_D(z) = \max(z)\;\mathrm{if}\;(z)\;\mathrm{terminal}$.  
This convention and
$\hat{t}_D(z) = \max(z)$ if $(z)$  terminal
are used to uniquely distinguish terminal vertices when the graphs are downward (definition~\ref{def:dwndrctgrph}).

We next give some standard graph theory terminology.
\begin{defn}
\label{def:termnot}
Let $z\in D$, $G_D = (D, \Theta_D)$.
The {\em adjacent vertices} of $z$, $G_D^z$, are defined by  $G_D^z= \{x\mid (z,x)\in \Theta_D\}.$
{\em Non-terminal adjacent vertices} of $z$, $NT_z$, are defined by 
$NT_z = \{x\mid x\in G_D^z,\, (x) \mathrm{\;not\;terminal}\}.$
The {\em terminal adjacent vertices} of $z$, $T_z$,  are defined by
$T_z=G_D^z \setminus NT_z$.
\end{defn}

\begin{defn}[\bfseries Recursive definitions $\hat{p}_D$, $\hat{t}_D$]
\label{def:recursivedefs}
Let $z\in D$, $G_D = (D, \Theta_D)$.  We define
$\hat{t}_D(z)= \max(z)$ if $(z)$ is terminal.
Else $\hat{t}_D(z)$ is the minimum over the set
\[
\{\hat{t}_D(x)\mid x\in NT_z\}\cup\{\min(x)\mid x\in T_z\}.
\]
We define $\hat{p}_D(z) = \max(z)$ if $(z)$ is terminal.
Else $\hat{p}_D(z)$ is the minimum over the set
\[
\{\hat{p}_D(x)\mid x\in NT_z\}\cup\{\min(x)\mid x\in T_z\}\cup\{\min(z)\}.
\]

{\bf Note:}  $\hat{t}_D(z)=\max(z)$  iff $(z)$ is terminal, and
$\hat{p}_D(z)=\max(z)$  iff $(z)$ is terminal.
\end{defn}

\begin{figure}[h]
\begin{center}
\includegraphics[scale=.85]{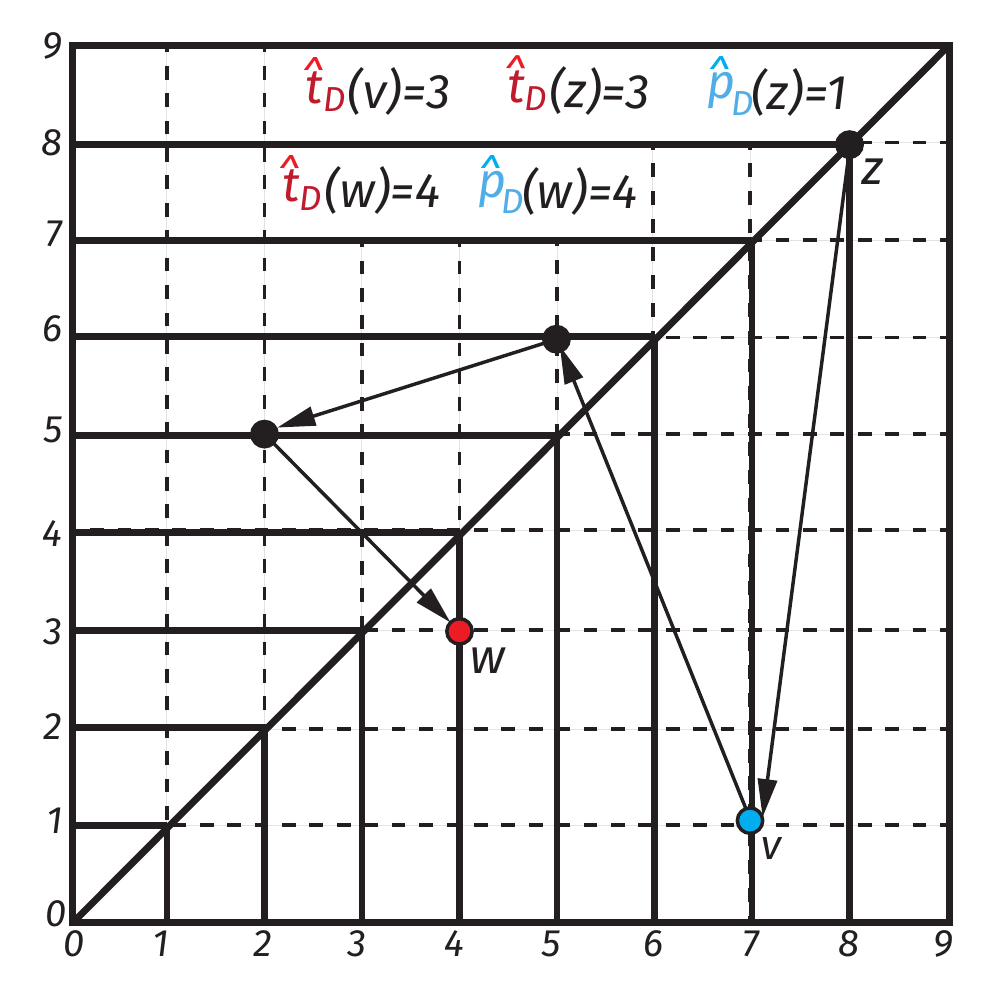}
\caption{Terminal vs. path labels, $k=2$}
\label{fig:termvspath}
\end{center}
\end{figure}

{\em Lattice Exit Story (some intuition):} 
We consider $G_D$ with terminal label function 
$\hat{t}_D$ (case $k=3$). Think of the digraph $G_D$ as a complex of caves (vertices) and tunnels (downward directed).  Each cave has a label on its wall giving its coordinates in $N^3$.  The value of  $\hat{t}_D(z)$ is also written on the wall of cave $z$.  An explorer is lowered into cave $z$ and tasked with finding a {\em terminal} cave  $x$, accessible from $z$, that has shortest distance to the boundary of the lattice (equal to $\hat{t}_D(z)$ and called the "lattice exit distance").  
Specifically, if $\hat{t}_D(z)\geq \min(z)$ then no exploration is needed; cave $z$ is closest to the boundary.  
If $\hat{t}_D(z)< \min(z)$ then the explorer finds a terminal path $z, \ldots, x$ with $\hat{t}_D(z)=\min(x)$. The path $(z, \ldots, x)$ is then written on the wall of cave $z$.  
The set $\{z\mid \hat{t}_D(z)< \min(z)\}$ gives the caves where getting out of $z$ gets one closer to the boundary.   The set $\{\hat{t}_D(z)\mid \hat{t}_D(z)<\min(z)\}$  represents the distances obtained by such explorations (i.e., {\em significant labels} for $\hat{t}_D$).

\begin{figure}[h]
\begin{center}
\includegraphics[scale=.85]{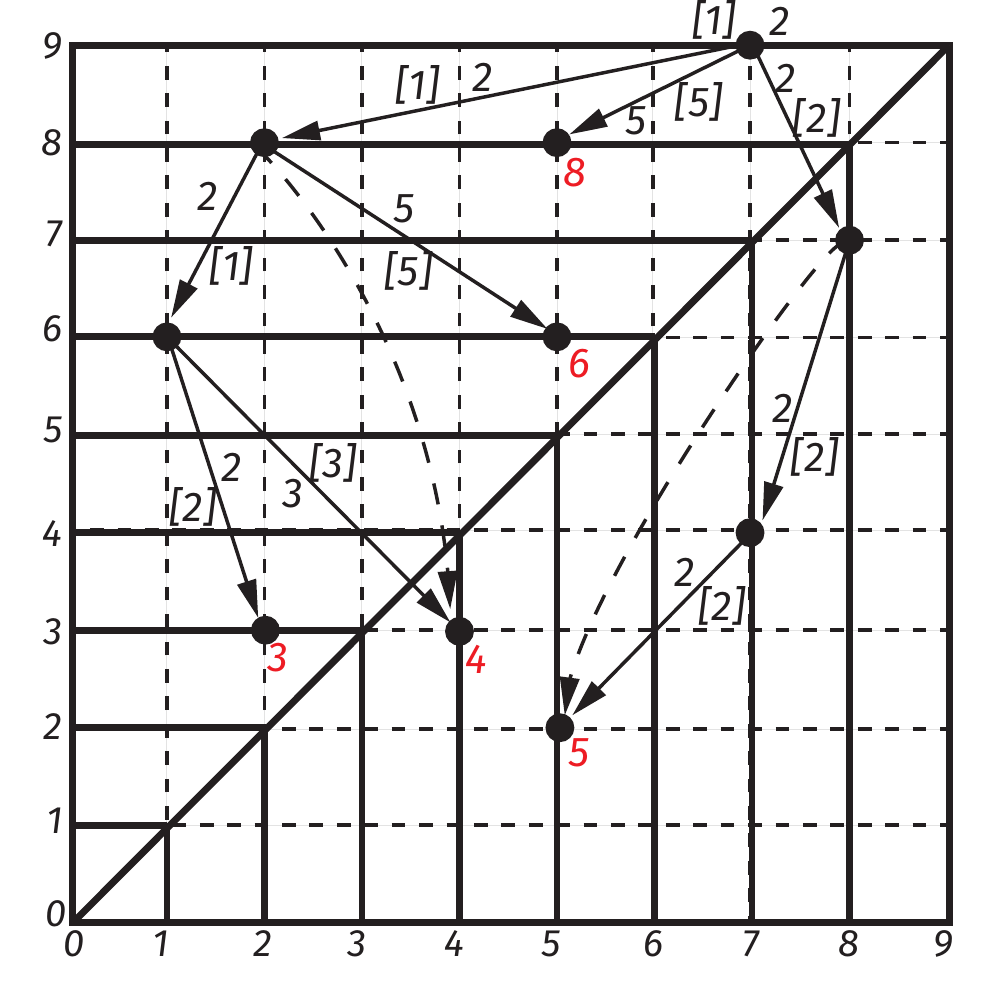}
\caption{${\bf [}\hat{p}_D(z){\bf ]}$ (in brackets) and $\hat{t}_D(z)$ by depth first search}
\label{fig:termvspath}
\end{center}
\end{figure}

\section{Regularities comfortably in ZFC}
We discuss the large scale regularties of $\hat{p}_D(z)$.  We start with a version of Ramsey's Theorem (\cite{rg:rt}, p.23).  First we define order equivalence of $k$-tuples.

\begin{defn}[\bfseries Equivalent ordered $k$-tuples]
\label{def:ordtypeqv}
Two k-tuples in $N^k$, $x=(n_1,\ldots,n_k)$ and $y=(m_1,\ldots,m_k)$, are  
{\em order equivalent tuples $(ot)$} if 
$\{(i,j)\mid n_i < n_j\} =  \{(i,j)\mid m_i < m_j\}$ and  $\{(i,j)\mid n_i = n_j\} =  \{(i,j)\mid m_i =m_j\}.$  
\end{defn}
Note that $ot$ is  an equivalence relation on $N^k$.
The standard SDR (system of distinct representatives) for the $ot$  
equivalence relation is gotten by replacing $x=(n_1,\ldots,n_k)$ by 
$\rho(x):= (\rho_{S_x}(n_1),\ldots, \rho_{S_x}(n_k))$ where
 $\rho_{S_x}(n_j)$ is the rank of $n_j$ in $S_x = \{n_1, \ldots, n_k\}$ (e.g,
$x=(3, 8, 5, 3, 8)$, $S_x = \{3, 5, 8\}$, $\rho(x)=(0, 2, 1, 0, 2)$).
The number of equivalence classes is $\sum_{j=1}^k \sigma(k, j)$ where
$\sigma(k,j)$ is the number of surjections from a $k$ set to  a $j$ set.

\begin{thm}[\bfseries Ramsey's theorem version]
If $f: N^r \rightarrow X$, 
$\rm{Im}(f)=\{f(z)\mid z\in N^r\}$ finite, then
there exists infinite $H=\{h_0, h_1, \ldots \}\subseteq N$ s.t. $f$ is constant on the order equivalence classes of $H^r$. 
\label{thm:ramsey}
\end{thm}

\begin{figure}[h]
\begin{center}
\includegraphics[scale=.75]{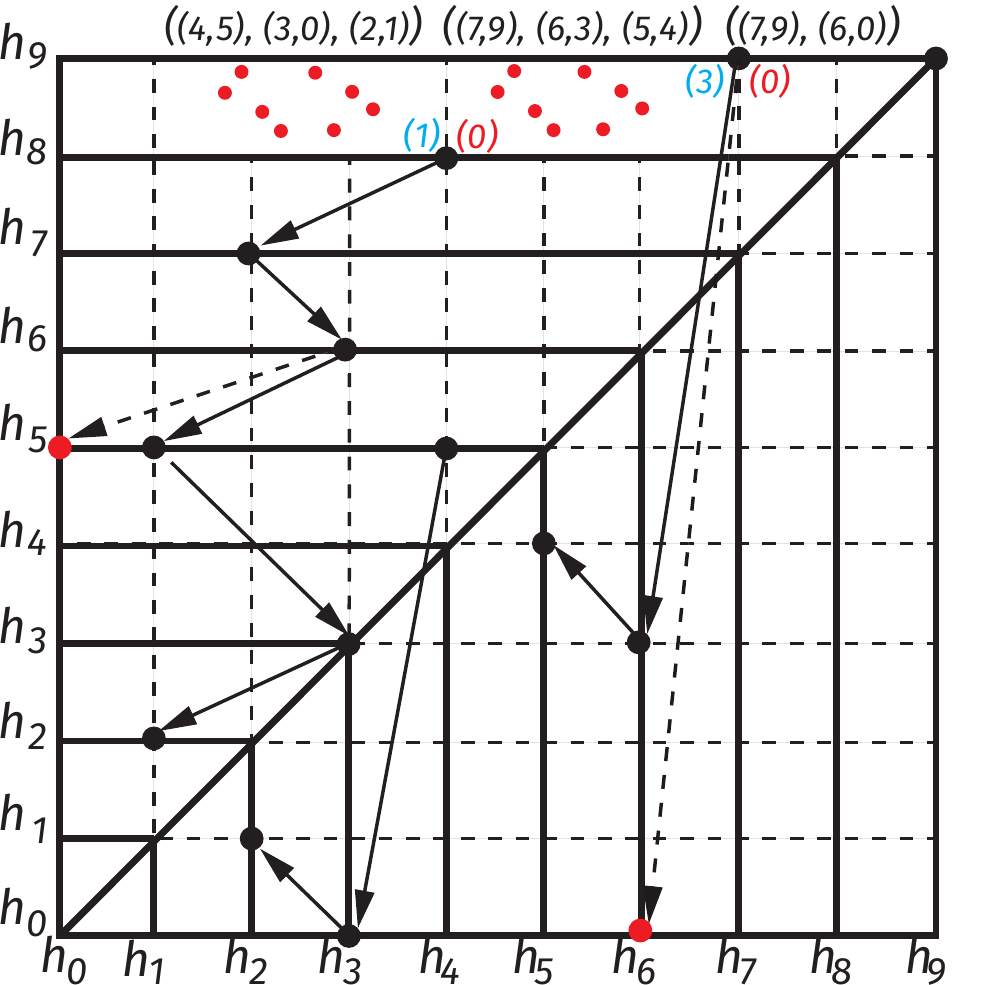}
\caption{Basic idea for proof of theorem~\ref{thm:regforhatp}}
\label{fig:hland}
\end{center}
\end{figure}

\begin{defn}[$\mathcal{X}$ notation]
Define $\mathcal{X}(statement) = 0\;{\rm if}\,statement\, {\rm false}$,
and $\mathcal{X}(statement) = 1\;{\rm if}\,statement\, {\rm true}$.
\end{defn}

\begin{thm}[\bfseries $\hat{p}_D$ large scale regularity structure]
\label{thm:regforhatp}
Let $G=(N^k, \Theta).$  
There exists an infinite $H=\{h_0, h_1, \ldots \}\subseteq N$ such that for all 
$z\in D= H^k$  either 
$\hat{p}_D(z) = h_0 < \min(z)$ or $\hat{p}_D(z)\geq \min(z).$ Thus, $\hat{p}_D	$ has at most one significant label:
$|\{\hat{p}_D(z)\mid \hat{p}_D(z)< \min(z)\}| \leq 1.$
\begin{proof}

Recall definition~\ref{def:pathlabel} of the function $\hat{p}_D(z)$.
Use Ramsey (\ref{thm:ramsey} with $r=2k$) on $N^{2k}$: $w=(n_1,\ldots  n_k, m_1,\ldots m_k)\in N^{2k}$ $\mapsto$
$(x, y)\in N^k \times N^k$ , $x=(n_1,\ldots n_k)$, $y=(m_1,\ldots m_k).$
Let $f(w)=\mathcal{X}((x,y)\in \Theta),$ so that $\mathrm{Im}(f)=\{0,1\}.$
By Ramsey's theorem, there is an infinite $H\subseteq N$ such that $f$ is 
constant on the order type equivalence classes of $H^{2k} \equiv H^k\times H^k$.
We show that $\hat{p}_D(z)$ has at  most one significant label on $D=H^k$.

Let $z=x_1, \ldots, x_t $ be a shortest path in $G_D$ from $z$ to a vertex 
$\min(x_t) = \hat{p}_D(z)$.   We claim $\min(x_t)=h_0$.  Otherwise, replace every minimum coordinate of $x_t$ by $h_0$ to obtain $\hat{x}_t$ and note that $(x_{t-1},x_t)$   and
$(x_{t-1},\hat{x}_t)$ have the same ot in $H^{2k}$.  
Thus, $(x_{t-1},\hat{x}_t)\in \Theta$ and $z=x_1, \ldots, \hat{x}_t $ is a shortest path
required.

\end{proof}
\end{thm}

\section{Basic definitions and theorems}

\begin{defn}[\bf decreasing sets of functions]
Let $f$ and $g$ be functions with domains contained in $N^k$ and ranges in $N$.
Define $f\geq g$ by 
$$(1)\;{\rm domain}(f)\subseteq {\rm domain}(g)\; and\;(2)\; {\rm for\; all}\;x\in {\rm domain}(f),\;f(x)\geq g(x).$$
A set $S$ of such functions is {\em decreasing} if for all $f,g\in S$ with
${\rm domain}(f)\subseteq {\rm domain}(g)$, $f\geq g$.
\end{defn}

\begin{defn}[\bf regressive value]
Let $X\subseteq N^k$ and $f:X\rightarrow Y\subseteq N$.  An integer $n$
is a  {\em regressive value} of $f$ on $X$
if there exist $x$ such that $f(x)=n<\min(x)$ .
\end{defn}

\begin{defn}[\bf field of a function and reflexive functions]
For $A\subseteq N^k$ define ${\rm field}(A)$ to be the set of all coordinates of elements of $A$.  A function $f$ is reflexive in $N^k$ if 
${\rm domain}(f) \subseteq N^k$ and  ${\rm range}(f) \subseteq {\rm field}({\rm domain}(f))$.
\end{defn}

\begin{defn}[the set of functions $T(k)$ ]
$T(k)$ denotes all reflexive functions with finite domain: $|{\rm domain}(f)|<\infty$.
\end{defn}

\begin{defn} [\bf full and jump free] 
Let $Q\subset T(k)$ denote a collection of reflexive functions in $N^k$ whose domains are finite subsets of $N^k$. 

\begin{enumerate}
\item {\bf full:}  $Q$ is a {\em full} family of functions on $N^k$ if for every finite subset 
$D\subset N^k$ there is at least one function $f$ in $Q$ whose domain is $D$.

 
\item{\bf jump free:} For $D\subset N^k$ and $x\in D$ define $D_x = \{z\mid z\in D,\, \max(z) < \max(x)\}$. 
Suppose that for all $f_A$ and $f_B$  in $Q$, where $f_A$ has domain $A$ and $f_B$ has domain $B$,  the conditions
 $x\in A\cap B$, $A_x \subseteq B_x$, and $f_A(y) = f_B(y)$ for all $y\in A_x$ imply that 
$f_A(x) \geq f_B(x)$.  Then $Q$ will be called a {\em jump free} family of functions on $N^k$.  
\end{enumerate}
\label{def:fullrefjf}
\end{defn}

\begin{defn}[\bf Regressively regular over $E$]
\label{def:regreg}
Let $k\geq 2$, $D\subset N^k$, $D$ finite, $f: D\rightarrow N$. 
We say $f$ is {\em regressively regular} over 
$E$, $E^k\subset D$, if for each $ot$ either (1) or (2):
\begin{enumerate}
\item{\bf decreasing mins:}  For all $x, y\,\in E^k$ of order type $ot$, $f(x)=f(y)< \min(E)$ 
\item{\bf non decreasing mins:}   For all $x\in E^k$ of order type $ot$ $f(x)\geq \min(x).$
\end{enumerate}
\end{defn}   

\begin{thm}[\bf Decreasing class]
\label{thm:decclass}
Let $k, p \geq 2$ and $S\subseteq T(k)$ be a full and decreasing family of functions.
Then some $f\in S$ has at most $k^k$ regressive values on some 
Cartesian power $E^k \subseteq {\rm domain}(f)$, $|E| = p$.  
In fact, there exists $E\subseteq A\subset N$, $|E|=p$ and $f\in S$,
${\rm domain}(f)=A^k$ such that
$f$ is regressively regular over $E$.
\end{thm}

\begin{thm}[\bf Jump free theorem (\cite{hf:alc}, \cite{hf:nlc})] 
\label{thm:jumpfree}
Let $p, k\geq 2$ and $S\subseteq T(k)$ be  a full and jump free family.
Then some $f\in S$ has at most $k^k$ regressive values on some 
$E^k \subseteq {\rm domain}(f)$, $|E| = p$.  
In fact, some $f\in S$ is regressively regular over some $E$ of cardinality $p$.
\end{thm}

\begin{figure}[h]
\begin{center}
\includegraphics[scale=.85]{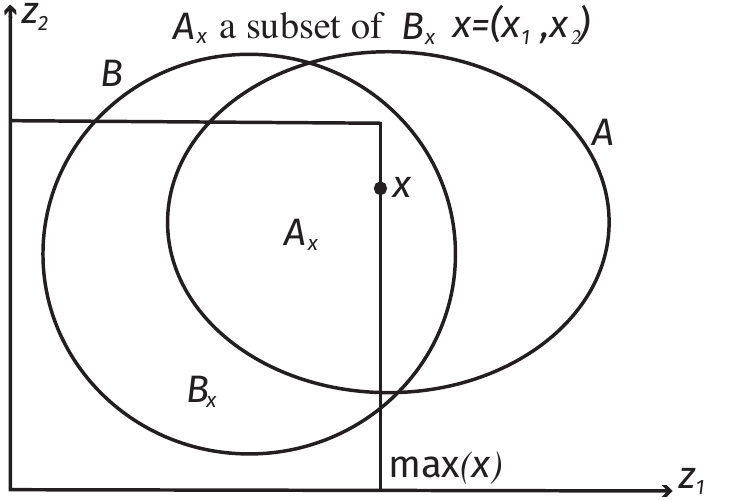}
\caption{Basic jump free condition~\ref{def:fullrefjf}}
\label{fig:jfvenn}
\end{center}
\end{figure}
 
\begin{figure}[h]
\begin{center}
\includegraphics[scale=.85]{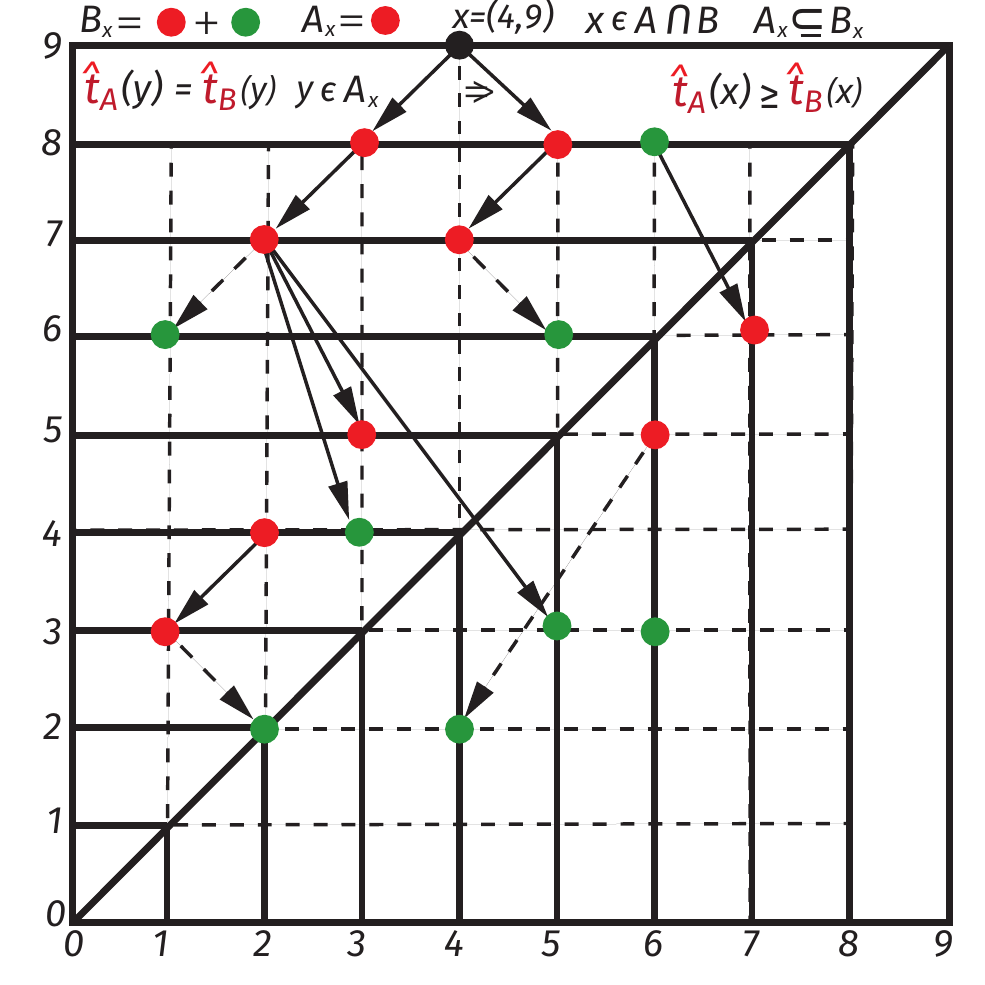}
\caption{Dashed edges not allowed by jump free~\ref{def:fullrefjf}}
\label{fig:jfexmp}
\end{center}
\end{figure}
We use ZFC for the axioms of set theory, Zermelo-Frankel  plus the axiom of choice (see Wikipedia).  
 The jump free theorem can be proved in 
 ZFC + ($\forall n$)($\exists$ $n$-subtle cardinal)
 but not in ($\exists$ $n$-subtle cardinal) for any fixed $n$ (assuming this theory is consistent).
A proof is in Section 2 of \cite{hf:alc},
``Applications of Large Cardinals to Graph Theory,'' October 23, 1997, No. 11 of  
\href{https://u.osu.edu/friedman.8/foundational-adventures/downloadable-manuscripts/}{Downloadable Manuscripts}.
The decreasing class theorem is proved in Section 1 of \cite{hf:alc} using techniques within ZFC (Ramsey theory in particular).
The jump free theorem is used to study lattice posets in \cite{jg:pos}
(Appendix A defines n-subtle cardinals).
The functions $\hat{p}_D$ define a full and decreasing  class and thus have large scale regularities of the form specified in 
Theorem~\ref {thm:decclass}.  
The functions $\hat{t}_D$ form a full but not decreasing class.
We will use the jump-free theorem to describe the large scale regularities of these functions.  

\section{Large scale regularities of more complex lattice exit models}

\begin{lem}[\bfseries $\{\hat{t}_D\}$ full, reflexive, jump free]
\label{lem:hattfrjf}
Take 
$$S=\{\hat{t}_D\mid D\subset N^k, |D|<\infty\}$$ (see~\ref{def:termlabel}).  
Then $S$ is full, reflexive, and jump free.

\begin{proof} See figures~\ref{fig:jfvenn} and~\ref{fig:jfexmp}.
Full and reflexive is immediate.
Let $\hat{t}_A$ and $\hat{t}_B$  satisfy the conditions of $f_A$ and $f_B$
in definition~\ref{def:fullrefjf}. 
Note that by definition, $x\notin A_x \,{\rm or}\; B_x$.
If $(x)$ is terminal in $A$ then $\hat{t}_A(x) =\max(x)\geq \hat{t}_B(x)$ by the downward condition on $G$.
Else, let $(x, \ldots, y)$ be a terminal path in $G_A$.  
Then $\hat{t}_B(y)=\hat{t}_A(y)=\max(y)$ implies  $(x, \ldots, y)$ is a
terminal path in $G_B$. Thus, $\hat{t}_A(x) \geq \hat{t}_B(x)$ as was to be shown.
\end{proof}
\end{lem}

\begin{figure}[h]
\begin{center}
\includegraphics[scale=.75]{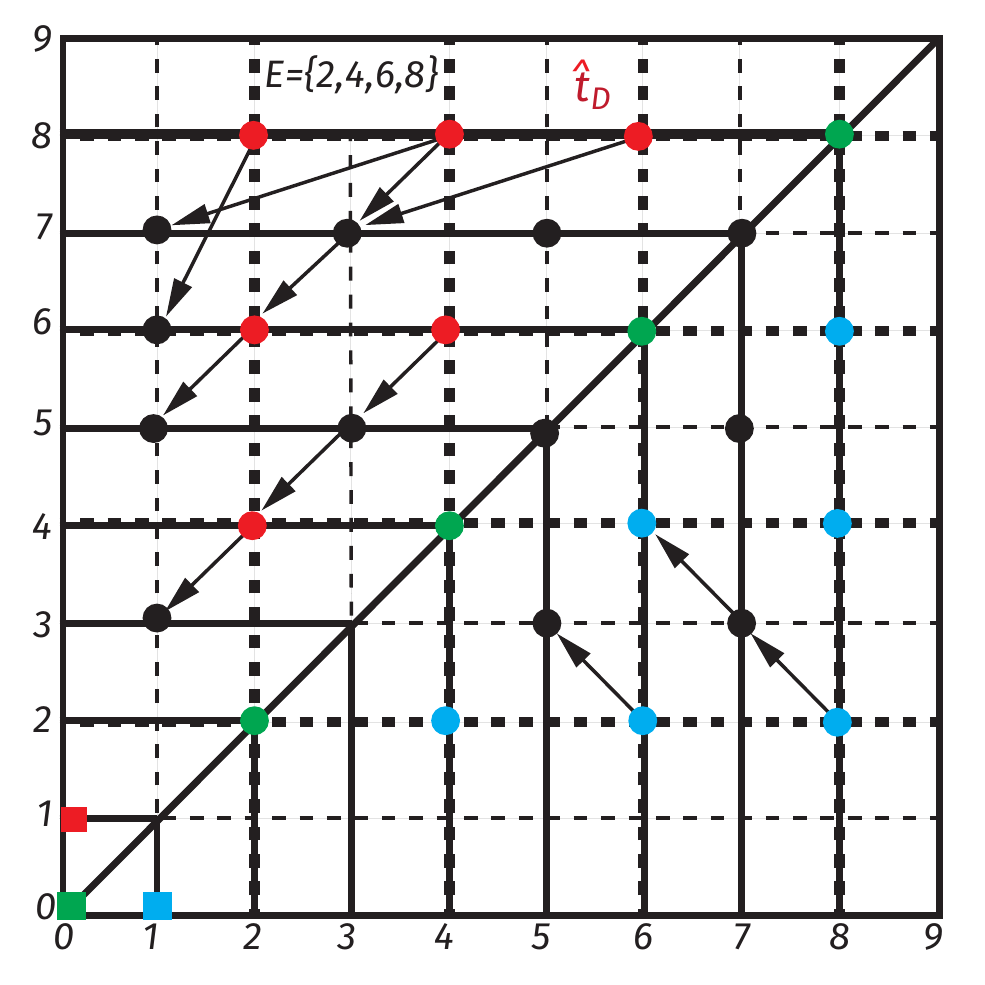}
\caption{Regressive regularity of $\hat{t}_D$ on $E^2=\{2,4,6,8\}^2$ \label{fig:rregtd}}
\end{center}
\end{figure}

\begin{thm}[\bfseries Jump free theorem for $\hat{t}_D$]
\label{thm:jfthat}
Let
$S=\{\hat{t}_D\mid D\subset N^k,\;|D|<\infty\}$ and
let $p, k\geq 2$.   Then some $f\in S$ has at most $k^k$ regressive values on some $E^k \subseteq {\rm domain}(f)$, $|E|=p$.  In fact, some $f\in S$ is regressively  regular over some $E$ of cardinality $p$.

\begin{proof}
Follows from lemma~\ref{lem:hattfrjf} and the jump free theorem~\ref{thm:jumpfree}.  
See figure~\ref{fig:rregtd} for an example of regressive regularity.
\end{proof}
\end{thm}


%
Following Friedman~\cite{hf:alc}:
\begin{defn} [{\bf Partial selection}]
\label{def:partselect}
A function $F$ with domain a subset of $X$ and range a subset of $Y$ will be called a {\em partial function}
from $X$ to $Y$ (denoted by $F: X\rightarrow Y$).  If $z\in X$ but $z$ is not in the domain of $F$, we say 
$F$ is {\em not defined} at $z$.
Let $r
\geq 1$.  A partial function 
$F: N^k\times(N^k \times N)^r \rightarrow N$
will be called a {\em partial selection} function if whenever 
$F(x, ((y_1,n_1), (y_2,n_2), \ldots (y_r,n_r)))$ is defined we have 
$F(x, ((y_1,n_1), (y_2,n_2), \ldots (y_r,n_r))) = n_i$ for some $1\leq i \leq r$.
\end{defn}

Next we generalize the $\hat{t}_D$ to a function $\hat{s}_D$.  We refer to the former function as the ''terminal vertex model'' and to the latter as the ''committee model." 

\begin{defn} [{\bf $\hat{s}_D$ for $G_D$}]
\label{def:chanlabel}
Let $r\geq 1,$ $z\in D$, $G_D = (D, \Theta_D)$, $D$ finite, $G_D^z= \{x\mid (z,x)\in \Theta_D\}.$ 
Let $F: N^k\times(N^k \times N)^r \rightarrow N$ be a partial selection function.  
We define $\hat{s}_D(z)$ recursively as follows. Let
\[
\Phi^D_z := \{ F[z, (y_1,n_1), (y_2,n_2), \ldots, (y_r,n_r)],\;y_i \in G^z_D\}
\]
be the set of defined values of $F$  where  
$n_i=\hat{s}_D(y_i)$ if $\Phi^D_{y_i}\neq\emptyset$ and
$\;n_i=\min(y_i)$ if $\Phi^D_{y_i}=\emptyset.\;$
If $\Phi^D_z=\emptyset$,  define $\hat{s}_D(z) = \max(z)$.
If $\Phi^D_z\neq\emptyset$,  define $\hat{s}_D(z)$ be the minimum over $\Phi^D_z$.
\end{defn}
 NOTE: If $\Phi^D_z\neq \emptyset$ then an induction on $\max(z)$ shows 
 $\hat{s}_D(z) < \max(z).$ Recall that $(G,\Theta)$ is downward.
Thus, $\Phi^D_z=\emptyset$ iff $\hat{s}_D(z) = \max(z)$.

\begin{thm}[\bfseries Large scale regularities for $\hat {s}_D$]
\label{thm:jfhats}
Let $r\geq 1$, $p, k\geq 2$.
$S=\{\hat {s}_D\mid D\subset N^k,\;|D|<\infty\}$. Then some 
$f\in S$ has at most $k^k$ regressive values over some 
$E^k \subseteq {\rm domain}(f)$, $|E|=p.$ 
In fact, some $f\in S$ is regressively  regular over some $E$ of cardinality $p$.

\begin{proof}
Recall \ref{thm:jumpfree}. Let  $S=\{\hat {s}_D\mid D\subset N^k,\;|D|<\infty\}.$
$S$ is obviously full and reflexive. 
We show $S$ is jump free.
We show for all $\hat {s}_A$ and $\hat {s}_B$  in $S$,  the conditions
 $x\in A\cap B$, $A_x \subseteq B_x$, and $\hat {s}_A(y) =\hat {s}_B(y)$ for all $y\in A_x$ imply that 
$\hat {s}_A(x) \geq \hat {s}_B(x)$.  (i.e., $S$ is  {\em jump free}). 
If $\Phi^A_x = \emptyset$ then $\hat {s}_A(x)=\max(x)\geq \hat {s}_B(x).$
Let $n=F[x, (y_1,n_1), (y_2,n_2), \ldots (y_r,n_r)]\in \Phi^A_x$ where 
$n_i=\hat {s}_A(y_i)$ if $\hat {s}_A(y_i)<\max(y_i)$ and 
$n_i=\min(y_i)$ if $\hat {s}_A(y_i)=\max(y_i).$
But $\hat {s}_A(y_i) =\hat {s}_B(y_i)$, $i=1, \ldots, r,$ implies
$n\in \Phi^B_x$ and thus $\Phi^A_x \subseteq \Phi^B_x$ and
$
\hat {s}_A(x)=\min(\Phi^A_x) \geq \min(\Phi^B_x)=\hat {s}_B(x).
$
\end{proof}
\end{thm}

\begin{figure}[h]
\begin{center}
\includegraphics[scale=.7]{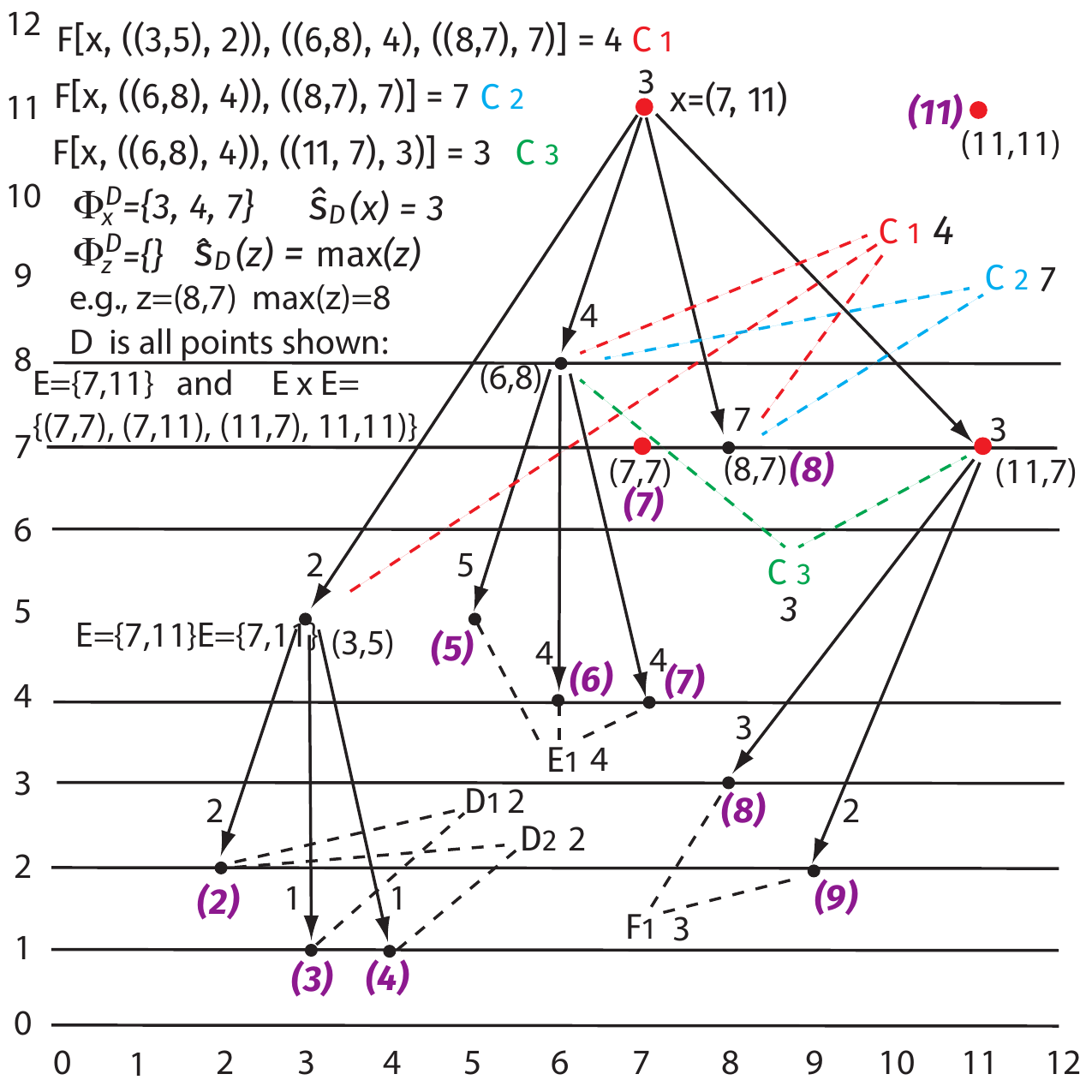}
\caption{An example of $\hat{s}_D$}
\label{fig:univs}
\end{center}
\end{figure}

As an example of computing $\hat{s}_D$, consider figure~\ref{fig:univs}.
The values of the terminal vertices where  $\Phi^A_x = \emptyset$ are shown in parentheses,  left to right: (2), (3), (4), (5), (6), (8), (8), (9).  These numbers are  $\max((a,b))$ for each terminal vertex $(a,b)$.
We assume we have a partial selection functions of the form 
$F:N^2 \times (N^2\times N)^r \rightarrow N$ ($r=2, 3$ here).  
To compute $\hat{s}_D(x)$ for $x=(7,11)$ there are three defined values:\\
$F[x, ((3,5),2), ((6,8),4), ((8,7),7)]=4$, \\
$F[x, ((6,8),4), ((8,7),7)]= 7$, \\
$F[x, ((6,8),4), ((11,7),3)]= 3$.\\ 
Intuitively, we think of these as (ordered) committees reporting values to the boss, $x=(7,11).$
The first committee, $\rm{C}1$, consists of subordinates, $(3,5),(6,8), (8,7)$
reporting respectively $2, 4, 7$.  
The committee decides to report $4$ (indicated by  $\rm{C}1\;4$ in  
figure~\ref{fig:univs}). The recursive construction starts with terminal vertices reporting their minimal coordinates.
But, the value reported by each committee is not, in general, the actual minimum of the reports of the individual members. 
Nevertheless, the boss, $x=(7,11)$  {\em always} takes the minimum of the values reported to him by the committees.  In this case the values reported by the committees are $4, 7, 3$ the boss takes $3$ (i.e., $\hat{s}_D(x)=3$ for the boss, $x=(7,11)$).

Observe in figure~\ref{fig:univs} that the values in parentheses, (2), (3), (4), (5), (6), (8), (8), (9), don't figure into the recursive construction of $\hat{s}_D$.  They immediately pass their minimum values on to the computation: 2, 1, 1, 5, 4, 4, 7, 3.
This leads to the following generalization of definition~\ref{def:chanlabel}.

\begin{defn} [{\bf $h^{\rho}_D$ for $G_D$}]
\label{def:genhat}
Let $r\geq 1$, $k\geq 2$, $z\in D$, $D$ finite, $G_D = (D, \Theta_D)$. 
Let $F: N^k\times(N^k \times N)^r \rightarrow N$ be a partial selection function.
Let $\rho_{D}: D \rightarrow N$ be such that $\min(x)\leq \rho_{D}(x)$.
We define $h^{\rho}_D$ recursively on $\max$. Let
\[
\Phi^D_z := \{ F[z, (y_1,n_1), (y_2,n_2), \ldots, (y_r,n_r)],\;y_i \in G^z_D\}
\]
be the set of defined values of $F$  where  
$n_i=h^{\rho}_D(y_i)$ if $\Phi^D_{y_i}\neq\emptyset,$ and
$\;n_i=\min(y_i)$ if $\Phi^D_{y_i}=\emptyset.\;$
If $\Phi^D_z\neq\emptyset$,  define $h^{\rho}_D(z)$ to be the minimum over 
$\Phi^D_z$.
If $\Phi^D_z=\emptyset$,  define $h^{\rho}_D(z) = \rho_{D}(z)$.
Note that $\rho_D$ need not be reflexive on $D$.
\end{defn}

\begin{lem}
\label{lem:shatvsh}
 Let $E$ be of cardinality $p$.
Then $\hat{s}_D$ is regressively  regular over $E$ iff $h^{\rho}_D$ is regressively  regular over $E$. In fact, $h^{\rho}_D(x)=\hat{s}_D(x)<\max(x)$ if $\Phi^D_x\neq \emptyset$.
\begin{proof}
Let $x, y\in E^k$.
From the recursive definitions~\ref{def:chanlabel}   and \ref{def:genhat}, the sets 
$\Phi^D_x$ are the same for  both 
 $h^{\rho}_D(x)$ and $\hat{s}_D(x)$. 
Thus,  $w=\hat{s}_D(x)=\hat{s}_D(y) < \min(E)$ iff $w=h^{\rho}_D(x)=h^{\rho}_D(y) < \min(E)$ as these relations imply both $\Phi^D_x\neq \emptyset$ and $\Phi^D_y\neq \emptyset$.
Likewise, if $\Phi^D_x\neq \emptyset$ then 
$\max(x)>\hat{s}_D(x)=h^{\rho}_D(x)\geq \min(x)$.
If $\Phi_x^D = \emptyset$, by definition $h^{\rho}_D(x) =\rho_{D}(x)\geq \min(x)$ and
$\hat{s}_D(x) = \max(x) \geq \min(x)$.
\end{proof}
\end{lem}

\begin{thm}[\bfseries Regressive regularity of $h^{\rho}_D$]
\label{thm:jfh}
Let $r\geq 1$, $p, k\geq 2$. Let $G=(N^k, \Theta)$ be downward directed.  Let $S=\{h^{\rho}_D\mid D\subset N^k,\;|D|<\infty\}$.  Then some $f\in S$ has at most $k^k$ regressive values on some $E^k \subseteq {\rm domain}(f)=D$, $|E|=p$.
In fact, some $f\in S$ is regressively  regular over some $E$ of cardinality~$p$.
\begin{proof}
Follows from lemma~\ref{lem:shatvsh} which shows that the sets, $E$, $|E|=p$, over which $h_D^\rho$ is regressively regular don't depend on the function $\rho_D$ as defined.   
In fact, $h^{\rho}_D(x)=\hat{s}_D(x)<\max(x)$ if $\Phi^D_x\neq \emptyset$.
If $\Phi^D_x= \emptyset$ then  the values $h^{\rho}_D(x)=\rho_D(x)$ are only constrained by the condition
$\rho_D(x)\geq \min(x)$.
\end{proof}
\end{thm}

Theorem~\ref{thm:jfhats}  with  $\max(x)$ replaced by $\min(x)$ when  $\Phi_x^D = \emptyset$  has been shown by Friedman to be independent of ZFC (same large cardinals as the jump free theorem).
See Theorem~4.4 through Theorem~4.15 \cite{hf:alc}.
Thus, a special case of theorem~\ref{thm:jfh} ($\rho_D=\min$) is independent of ZFC.  
Lemma~\ref{lem:shatvsh} shows that theorem~\ref{thm:jfh} for any $h^{\rho}_D$ results in exactly the same sets $E$ of regressive regularity as theorem~\ref{thm:jfhats}.
Hence,
theorem~\ref{thm:jfh} provides a family of ZFC independent jump free type theorems parameterized by  the $\rho_D$.

\begin{figure}[h]
\begin{center}
\includegraphics[scale=0.93]{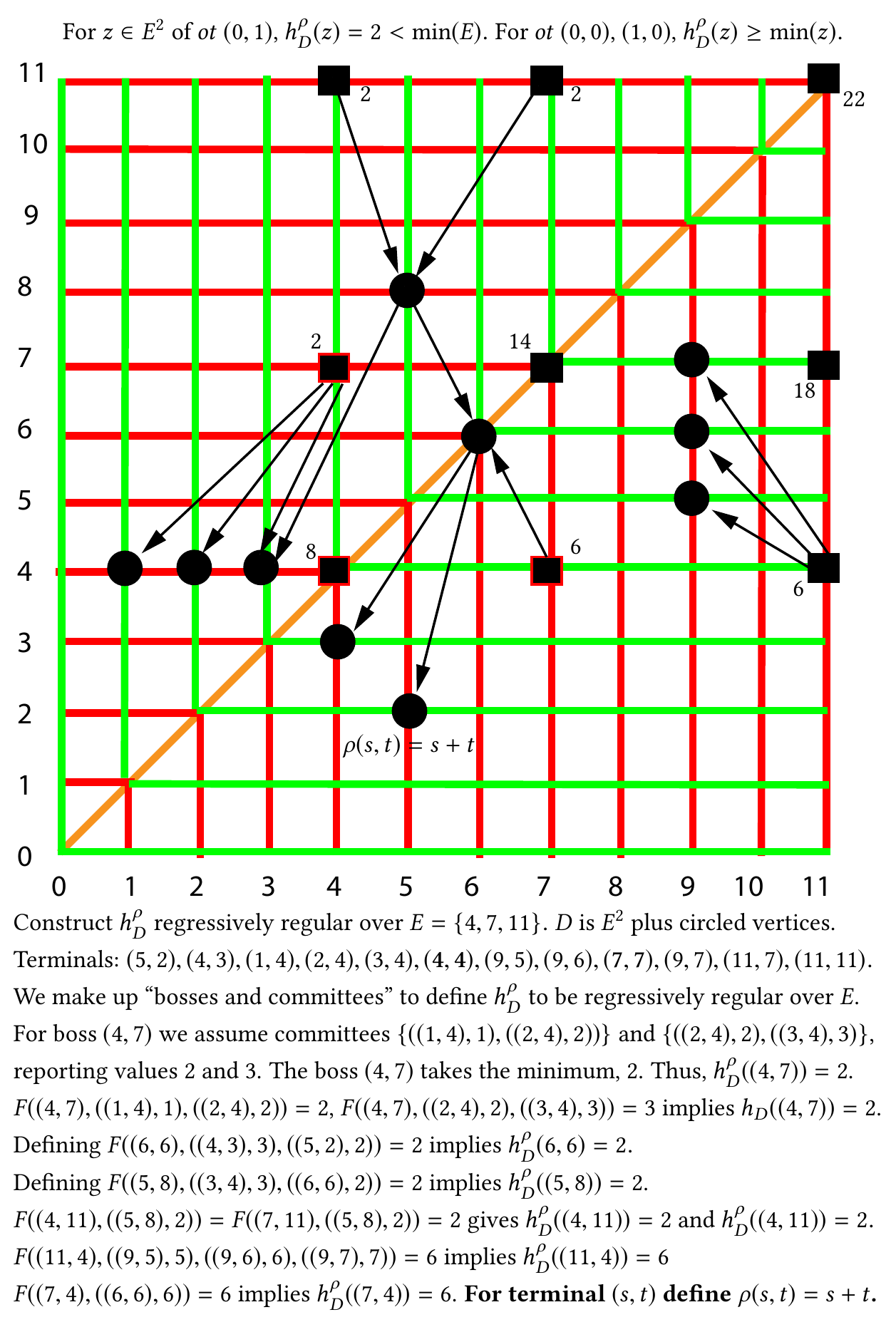}
\caption{Computing example of $h_D^{\rho}$}
\label{fig:exhdrho}
\end{center}
\end{figure}

Figure~\ref{fig:exhdrho} shows an example of a regressively regular $h_D^{\rho}$ function over a set $E$.  The choices of $F$ (definition~\ref{def:partselect})  aren't unique.
We use $\rho\coloneq\rho_D$.
For $z=(x,y)$ we define $\rho(z)= x+y$, labeling only $E^2$ with these.

\section{Combinatorial formulations}
\begin{figure}[h]
\begin{center}
\includegraphics[scale=0.95]{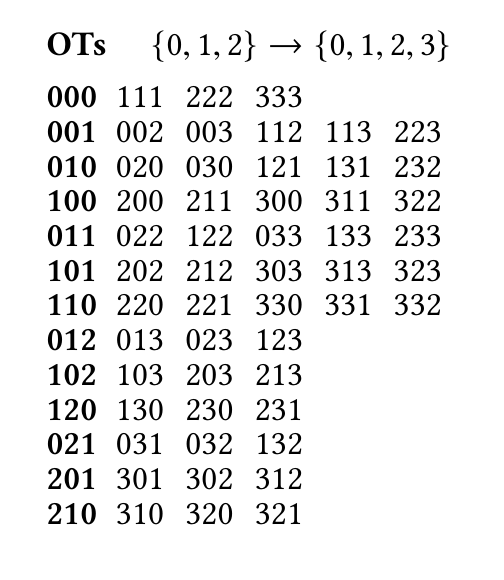}
\caption{Canonical order type array $T$ entries in $\bF_{3,4}$}
\label{fig:otypes}
\end{center}
\end{figure}

We start by defining some canonical systems of distinct representatives (SDRs).
\begin{defn}[SDRs]
\label{def:sdrs}
Let ${\rm SUR}(k,j)$ be the surjective maps $\{0, \ldots, k-1\}$ to
$\{0, \ldots, j-1\}.$  Let $\bF_{kp}=\{f\mid f:\{0, \dots k-1\}\rightarrow \{0,\ldots, p-1\}\}.$ 
Define
$${\rm OT(k,p)}\coloneq\cup_{j=1}^{p} {\rm SUR}(k,j).$$
${\rm OT(k,p)}$ is the {\em canonical SDR}  for the order type equivalence relation on $\bF_{kp}$.  
Let $E=\{e_0, \ldots, e_{p-1}\}$ be a subset of $N$, $e_0< \cdots <e_{p-1}.$
Let 
 $E^k=\{e_f\mid  f\in \bF_{kp}\}$  where $e_f= (e_{f(0)}, \dots, e_{f(p-1)})$.
Define
$${\rm OT(k,p,E)}= \{e_{f}\mid f\in {\rm OT(k,p)}\}$$ 
to be the canonical SDR for order type equivalence on $E^k$.
Define $\bE^k_{f}= \{e_g \mid g\in \bF_{kp},\; g \sim f\}$ to be the equivalence class in 
$E^k$ associated with $f\in {\rm OT(k,p)}$. 

Referring to figure~\ref{fig:otypes}, 
$\bE^3_{120} = \{(e_1,e_2,e_0), (e_1,e_3,e_0),(e_2,e_3e_0),(e_2,e_3,e_1)\}.$
\end{defn}

\begin{defn}[ Canonical order type array $T$] 
\label{def:canonical}
Define a canonical {\em order type array}, $T$, with entries in $\bF_{kp},$ as follows.
The first column of $T$ is the vector 
$T^{(1)}=(f_{11}, f_{21}, \ldots f_{m1})$ were the $f_{i1}$ are the elements of 
$\mathrm{OT}(k,p).$ The sets ${\rm SUR}(k,j)$ are listed by $j$ with elements of each set in lexicographic order. 
The row, $T_{(i)}$, $i=1, \ldots, m$, is 
$T_{(i)} = ( f_{i1}, f_{i2}, \dots, f_{in_i})$, the lexicographically ordered equivalence class associated with $f_{i1}$.
The entries $T(i,j)=f_{i j}$ where $i=\pi(f_{i1})$, $j=\pi(f_{ij})$, are the 
positions of $f_{i1}$ and $f_{ij}$ in the  lists of column $T^{(1)}$
and row $T_{(i)}$ respectively.
\end{defn}

\begin{defn}[$T_E^k$ and $gT_E^k$]
\label{def:arrays}
Let $T$ be a canonical order type array over $\bF_{kp}$, $E=\{e_0, \ldots, e_{p-1}\}\subset N$.
Replacing each $f$ in $T$ by $e_f$ (definition~\ref{def:sdrs}) gives an array which we denote by $X=T_E^k$. Thus, $X(i,j)=e_{f_{ij}}.$
If  $\mathrm{domain}(g)\supseteq E^k$ then $Y=gT_E^k$ has 
$Y(i,j)=g (e_{f_{ij}}) = g(X(i,j))$.  
\end{defn}

Consider figure~\ref{fig:exhdrho} where $E=\{4, 7, 11\}$. 
\[
T=\;
\begin{matrix}
00 & 11& 22\\
01 & 02 & 12\\
10 & 20 & 21
\end{matrix}
\hspace{.5cm}
T^2_{E}=\;
\begin{matrix}
e_{00} & e_{11}& e_{22}\\
e_{01} & e_{02} & e_{12}\\
e_{10} & e_{20} & e_{21}
\end{matrix}
\;=\;
\begin{matrix}
(4,4) & (7,7) & (11,11)\\
(4,7) & (4,11) & (7,11)\\
(7,4) & (11,4) & (11,7)
\end{matrix}
\]
Using $h_D^{\rho}$ as in figure~\ref{fig:exhdrho} we get
\[
h_D^\rho T^2_{E} =
\begin{matrix}
h_D^{\rho}(4,4) & h_D^{\rho}(7,7) & h_D^{\rho}(11,11)\\[0.4em]
h_D^{\rho}(4,7) & h_D^{\rho}(4,11) & h_D^{\rho}(7,11)\\[0.4em]
h_D^{\rho}(7,4) & h_D^{\rho}(11,4) & h_D^{\rho}(11,7)
\end{matrix}
\;=\;
\begin{matrix}
8 & 14 & 22\\
2 & 2 & 2\\
6 & 6 & 18
\end{matrix}
\]
which displays the regressive regularity over $E$.

Note that the pair of arrays, $(X,Y)$, where $X=T^2_{E}$ and $Y=h_D^\rho T^2_{E} ,$

\[
X=\;
\begin{matrix}
(4,4) & (7,7) & (11,11)\\
(4,7) & (4,11) & (7,11)\\
(7,4) & (11,4) & (11,7)
\end{matrix}
\;\;\;\mathrm{and}
\;\;\;
Y=\;
\begin{matrix}
8 & 14 & 22\\
2 & 2 & 2\\
6 & 6 & 18
\end{matrix}
\]
completely specifies the function $h_D^{\rho}$.

We use the notation of definition~\ref{def:genhat}, lemma~\ref{lem:shatvsh} and
theorem~\ref{thm:jfh}.

\begin{defn}[$D$  capped by $E^k\subset D$]
\label{def:cap}
For $D\subset N^k$, let $\max(D)$ be the maximum over $\max(z)$, $z\in D$.  
Let $\setmax(D)=\{z\mid z\in D, \max(z)=\max(D)\}$.
If  $\setmax(D) = \setmax(E^k)$, we say that $D$ is {\em capped by} $E^k\subset D$ with 
the {\em cap} defined to be $\setmax(E^k)$.
\end{defn}
 
\begin{defn}[Canonical capped bi-array]
\label{def:canordtyp}
Let $T$ be a canonical order type array over $\bF_{kp}$, $E=\{e_0, \ldots, e_{p-1}\}\subset N$.
The pair $(X,Y)=(T^k_E, h_D^{\rho}T^k_E)$ is a {\em canonical bi-array} representation of the function $h_D^{\rho}$.
Let $D$ be capped by $E^k\subseteq D$.
The pair $(X,Y)=(T^k_E, h_D^{\rho}T^k_E)$ is a {\em canonical capped bi-array} representation of the function $h_D^{\rho}$.
We write $(X,Y)_p=(T^k_E, h_D^{\rho}T^k_E)$ to indicate $|E|=p$.
\end{defn}

See figure~\ref{fig:rregtd} for an example where  $D$ is capped by $E^k\subseteq D$.
The downward condition on $G$ implies it is always possible to choose $D$ to satisfy this condition
without changing the function $(X,Y)=(T^k_E, h_D^{\rho}T^k_E)$.  Such a   $D$, capped by $E^k$, contains a description of $E^k$ exposed in the cap,
$\setmax(E^k).$

Note the following invariants of $(X,Y)=(T^k_E, h_D^{\rho}T^k_E)$. 
Any pair $(X_\sigma, Y_\sigma)$, where the rows are permuted by 
$\sigma$ also represents the function $h_D^{\rho}$ and preserves the elements of the first column.   
Given a permutation $\tau$ on $\{1, \ldots, n_i\}$  
with $\tau(1)=1$ the elements of rows $X_{(i)}$ and $Y_{(i)}$, $i=1, \ldots m_i$,  can be replaced by 
$X(i,\tau(1)), \ldots, X(i,\tau(i))$ and $Y(i,\tau(1)), \ldots, Y(i,\tau(i))$, preserving the elements of the first column while still representing the function $h_D^{\rho}$.

We use the notation of definition~\ref{def:genhat}, lemma~\ref{lem:shatvsh} and
theorem~\ref{thm:jfh}.
\begin{defn} [subsets of $E^k$]
\label{def:dompart} 
 Let $f$ be regressively regular over
$E\subset N$, $|E|=p$
(definition~\ref{def:regreg}). 
Define subsets
\[
E_L^k\coloneq\{x\mid x\in E^k, f(x) < \min(E)\},\;
E_U^k\coloneq\{x\mid x\in E^k, f(x)\geq \min(x)\}.
\] 
For $f= h_D^{\rho}$,
$E_L^k=\{x\mid x\in E^k, h_D^{\rho}(x) < \min(E)\}$.  $E_U^k$ is further partitioned
\[
E^k_{\neq\emptyset}=\{ x\mid x\in E_U^k, \Phi_D(x)\neq \emptyset \}\;\mathrm{and}\;
E_{\emptyset}^k=\{x\mid x\in E_U^k,\Phi_D(x) = \emptyset\}.
\]
\end{defn}

Note that both $E_L^k$ and $E_U^k$, when nonempty, are unions of elements (blocks) of ${\rm OT(k,p,E)}$ (definition~\ref{def:regreg}).
By definition, the regressively regular $f$ is constant on the blocks of the order type equivalence classes contained in $E_L^k$.  Recall that $\rho_D$ is restricted to 
$E_{\emptyset}^k$ in the recursive construction of $h_D^{\rho}$ and can be changed on this set as long as the condition $\rho_D(z) \geq \min(z)$ holds. Such changes in $h_D^{\rho}$ leave the sets of definition~\ref{def:dompart} invariant.

As an example, consider figure~\ref{fig:exhdrho}.
There, $k=2$, $p=3$ with $E=\{4,7,11\}$; $E^2$ is indicated by small squares, $D$ by squares plus circles.  
The SDR for order type equivalence on $\bF_{2,3}$ is \\
${\rm OT}(2,3) =\{(0,0), (0,1), (1,0)\}$, \\
${\rm OT}(2,3,E) =$
\[
 \{(e_{(0,0)}, e_{(0,1)}, e_{(1,0)}\}=
\{(e_0, e_0), (e_0, e_1),(e_1, e_0)\}=
\{(4, 4), (4, 7),(7, 4)\}.
\]

We have $E^k=\{e_f \mid f\in \bF_{2,3}\}$ (figure~\ref{fig:exhdrho}.)
$E^k_L=\{(4,7), (4,11), (7,11)\}$ and, in this case, the order equivalence class
$\bE^k_{(0,1)} = E^k_L$.  
In this example, $E^k_U = \bE^k_{(0,0)} \cup \bE^k_{(1,0)}$ with
 $$
 E^k_\emptyset = \mathrm{diag}(E^k)\cup \{ (11,7)\}\; \mathrm{and}\; 
 E^k_{\neq\emptyset} = E^k_L \cup  \{(11,4), (7,4)\}
 $$
where 
$\mathrm{diag}(E^k) = \{e_{00}, e_{11}, e_{22}\} = \{(4,4), (7,7),(11,11)\}.$
It is easy to see in general that $E^k_L \subseteq E^k_{\neq\emptyset}$ and that either $\mathrm{diag}(E^k) \subset E^k_{\emptyset}$ or
$\mathrm{diag}(E^k) \subset E^k_L$.

\section{Using the flexibility of $\rho_D$}
\begin{defn}[Regressive regularity for $(X,Y)$]
A canonical (capped) bi-array $(X,Y)=(T^k_E, h_D^{\rho}T^k_E)$ is {\em regressively regular} if 
for each row $Y_{(i)}$, $i=1,\ldots, m$, either
\[
 Y_{(i)}(j)\geq \min(X_{(i)}(j)),\;1\leq j \leq n_i\]
 or
\[ Y_{(i)}(s) = Y_{(i)}(t) < \min(E),\;1\leq s,t\leq n_i.
\]
\end{defn}

\begin{thm}[\bfseries Version of theorem~\ref{thm:jfh} for bi-arrays]
\label{thm:jfhc}
Let $r\geq 1$, $p, k\geq 2$. Let $G=(N^k, \Theta)$ be downward directed.
There is a canonically capped bi-array 
$(X,Y)=(T^k_E, h_D^{\rho}T^k_E)\,$, $|E|=p,\,$ such that 
$
|\{Y(i,j)\mid Y(i,j)< \min(X(i,j))\}|<k^k.
$
In fact, there exists a regressively regular canonical capped bi-array 
$(X,Y)=(T^k_E, h_D^{\rho}T^k_E)$, $|E|=p$. 
\begin{proof}
Restatement of theorem~\ref{thm:jfh} for  canonical capped bi-arrays.  
\end{proof}
\end{thm}


\begin{lem}[Choosing $\rho_D$]
\label{lem:capba}
Let $r\geq 1$, $p, k\geq 2$. Let $G=(N^k, \Theta)$ be downward directed.
Let $(X,Y)=(T^k_E, h_D^{\rho}T^k_E)$, $|E|=p$,
be a regressively regular canonical capped bi-array.  
Assume $E^k_L$ is nonempty and $\diag(E^k)\subseteq E^k_{\emptyset}$.
Let $e_{\vec{\mathbf{0}}}=(e_0, e_0, \ldots, e_0)\in E^k$.
Then $\rho_D$ can be chosen such that
\[
(1)\;\;\sum_{z\in E^k_L}h_D^{\rho}(z)+h_D^{\rho}(e_{\vec{\mathbf{0}}})
=(|E^k_L|+1)e_0.
\]
\begin{proof}
We can use the notation of either theorem~\ref{thm:jfhc} or theorem~\ref{thm:jfh}.
We use the latter.
Recall that $E^k_L = \{z\mid h_D^{\rho}(z)<e_0\}$ where $e_0=\min(E).$  Setting
$S=\sum_{z\in E^k_L}h_D^{\rho}(z)$ we have $S<|E^k_L|e_0$.
But $\diag(E^k)\subseteq E^k_{\emptyset}$ implies 
 $ h_D^\rho(e_{\vec{\mathbf{0}}})=\rho_D(e_{\vec{\mathbf{0}}})$
which can be assigned any value 
$\rho_D(e_{\vec{\mathbf{0}}})\geq \min(e_{\vec{\mathbf{0}}}) = e_0.$
Thus, assign $\rho_D(e_{\vec{\mathbf{0}}})=e_0+(|E^k_L|e_0-S)>e_0$ since 
$E^k_L\neq \emptyset$. Thus $(1)$ is satisfied.
\end{proof}
\end{lem}


As an aside, in the notation of theorem~\ref{thm:jfhc} we can calculate $|E^k_L|e_0$ as follows:
\[
(1)\;\;|E^k_L|e_0 =e _0 \sum_{i=1}^{m}  \mathcal{X}(Y(i,1)<e_0)n_i.
\]
Also, 
\[
(2)\;\; \sum_{z\in E^k_L}h_D^{\rho}(z) =\sum_{i=1}^{m}  \mathcal{X}(Y(i,1)<e_0)
 \sum_{j=1}^{n_i}Y(i,j).
\]


We engage in a ``thought experiment'' by using theorem~\ref{thm:jfhc} and 
lemma~\ref{lem:capba}  to construct a class of sequences of  instances to the classical subset sum problem.   
We assume
that for each $z\in \diag(N^k)$, the set of partial selection functions (see~\ref{def:chanlabel}) of the form 
$
F[z, (y_1,n_1), \ldots ], \,z \in \diag(N^k),\,
$
is empty.  
This {\em restricted diagonal} condition guarantees that, $r\geq 1$, 
\[
\Phi^D_z = \{ F[z, (y_1,n_1), (y_2,n_2), \ldots, (y_r,n_r)],\;y_i \in G^z_D\}
\]
is empty for $z\in \diag(E^k)$.
This implies $\diag(E^k)\subseteq E^k_{\emptyset}$ as in lemma~\ref{lem:capba}.

Let $r\geq 1$, $p, k\geq 2$. Let $G=(N^k, \Theta)$ be downward directed and diagonally restricted.
Let $\;(X,Y)=(T^k_E, h_D^{\rho}T^k_E)$, $|E|=p=2, 3, \ldots\,,$ be  a sequence of regressively regular canonical capped bi-arrays.

With each  bi-array we associate the multiset 
\[
\bM_p=\{Y(i,j)\mid 1\leq i \leq m, 1\leq j\leq n_i, Y(i,1)<e_0 \}
 \cup \{Y(1,1), \ldots Y(1,p)\}.
\]
Assume first that $E^k_L\neq \emptyset$. 
For each such $\bM_p$ let $t_p=(|E^k_L|+1)e_0$ as in lemma~\ref{lem:capba}.
As in lemma~\ref{lem:capba}, take 
$Y(1,1)=\rho_D(e_{\vec{\mathbf{0}}})$ 
 as specified in lemma~\ref{lem:capba}.
Lemma~~\ref{lem:capba} states that  for this instance there is a solution to the subset sum problem.
Choose $Y(1,2), \ldots, Y(1,p)$ such that this solution is unique.

If $E^k_L=\emptyset$, choose $Y(1,1), \ldots, Y(1,p)$ and $t_p$ such that there is no solution.

Thus, the instances to the subset sum problems just described have solutions if and only if $E^k_L\neq \emptyset$.  This condition can be verified by inspecting  the first column, $Y^{(1)}$ and comparing it with the first column $X^{(1)}$ (because of regressive regularity).  The number of comparisons is $m<k^k$, $k$ fixed.  

To summarize, we have defined a class of instances to the subset sum problem, parameterized by $r\geq 1$, $k\geq 2$ and $G=(N^k, \Theta)$, a  downward directed and diagonally restricted graph.  The parameter $p$ goes to infinity to measure the size of the instances. We fix $k$ and vary $r$ and $G$. 
Our procedure for checking the solutions for these instances is localized to the first columns of $X$ and $Y$ and thus bounded by $k^k$.

The existence of these instances and their solution has been demonstrated by using a corollary to a theorem independent of ZFC  (theorem~\ref{thm:jfhc}).
No other proof of existence is known to us.  

The corollary we used is  theorem~\ref{thm:jfhc} with  ``downward directed graph $G=(N^k, \Theta)$'' replaced by ``downward directed and diagonally restricted graph $G=(N^k, \Theta).$'' We conjecture that this corollary is also independent of ZFC.   In this case, the only proof of the existence of these instances and their solution would be from a ZFC independent theorem.

{\bf Acknowledgments:}  The author thanks  Professors Jeff Remmel and Sam Buss (University of California San Diego, Department of Mathematics) and 
Professor Emeritus Rod Canfield (University of Georgia, Department of Computer Science) for their helpful comments and suggestions.

%


\begin{thebibliography}{GRS90}

\bibitem[Fri97]{hf:alc}
Harvey Friedman.
\newblock Applications of large cardinals to graph theory.
\newblock Technical report, Department of Mathematics, Ohio State University,
  1997.

\bibitem[Fri98]{hf:nlc}
Harvey Friedman.
\newblock Finite functions and the necessary use of large cardinals.
\newblock {\em Ann. of Math.}, 148:803--893, 1998.

\bibitem[GRS90]{rg:rt}
R.~L. Graham, B.~L. Rothschild, and J.~H. Spencer.
\newblock {\em Ramsey Theory 2nd Ed.}
\newblock John Wiley, New York, 1990.

\bibitem[RW99]{jg:pos}
Jeffrey~B. Remmel and S.~Gill Williamson.
\newblock Large-scale regularities of lattice embeddings of posets.
\newblock {\em Order}, 16:245--260, 1999.

\end{thebibliography}

\end{document}